\documentclass[aos,nameyear,dvips]{arximspdf}

\doi{10.1214/10-AOS849}
\volume{39}
\issue{2}
\pubyear{2011}
\firstpage{1012}
\lastpage{1041}

\makeatletter

\newtheorem{theorem}{Theorem}
\newtheorem{corollary}{Corollary}
\newtheorem{proposition}{Proposition}
\newtheorem{lemma}{Lemma}

\newproclaim{remark}{Remark}

\newcommand{\RR}{\mathbb{R}}
\newcommand{\goto}{\rightarrow}
\newcommand{\hf}{{1 \over2}}

\newcommand{\V}{\operatorname{Var}}

\makeatother

\begin{document}
\begin{frontmatter}

\title{Testing composite hypotheses, Hermite polynomials and
optimal estimation of a nonsmooth functional}
\runtitle{Estimation of a nonsmooth functional}

\begin{aug}
\author[A]{\fnms{T. Tony} \snm{Cai}\thanksref{t2}\ead[label=e1]{tcai@wharton.upenn.edu}}
and
\author[A]{\fnms{Mark G.} \snm{Low}\corref{}\ead[label=e2]{lowm@wharton.upenn.edu}}
\runauthor{T. T. Cai and M. G. Low}
\affiliation{University of Pennsylvania}
\address[A]{Department of Statistics\\
The Wharton School\\
University of Pennsylvania\\
Philadelphia, Pennsylvania 19104\\
USA\\
\printead{e1}\\
\phantom{E-mail: }\printead*{e2}} 
\end{aug}

\thankstext{t2}{Supported in part by NSF FRG Grant DMS-08-54973.}

\received{\smonth{2} \syear{2010}}
\revised{\smonth{8} \syear{2010}}

%
\begin{abstract}
A general lower bound is developed for the minimax risk when estimating
an arbitrary functional. The bound is based on testing two composite
hypotheses\vspace*{1pt} and is shown to be effective in estimating the nonsmooth
functional ${1 \over n} \sum|\theta_i|$ from an observation $Y \sim
N(\theta, I_n)$. This problem exhibits some features that are
significantly different from those that occur in estimating
conventional smooth functionals. This is a setting where standard
techniques fail to yield sharp results.

A sharp minimax lower bound is established by applying the general
lower bound technique based on testing two composite hypotheses. A key
step is the construction of two special priors and bounding the
chi-square distance between two normal mixtures. An estimator is
constructed using approximation theory and Hermite polynomials and is
shown to be asymptotically sharp minimax when the means are bounded by
a given value $M$. It is shown that the minimax risk equals
$\beta_*^2M^2 ({\log\log n \over\log n})^2$ asymptotically, where
$\beta_*$ is the Bernstein constant.

The general techniques and results developed in the present paper can
also be used to solve other related problems.
\end{abstract}

%

\end{frontmatter}

\section{Introduction}\label{intro}

Minimax risk is one of the most commonly used benchmarks for
evaluating the performance of any estimation method. For this reason
considerable
effort has been made developing minimax theories in the
nonparametric function estimation literature. A key step in
all these developments is the derivation of minimax lower bounds.
Several effective lower bound techniques based on testing have been
introduced in the literature,
and it is often sufficient to derive the optimal
rate of convergence based on testing a pair of simple
hypotheses. Le Cam's method is a well-known approach based on this
idea. See, for example, \citet{LeCam73} and \citet{DL2}.

For estimation of quadratic functionals the story is somewhat more
complicated.
If the parameter space is not too
``large,'' regular parametric rate of convergence can be attained.
However \citet{BR88} showed that when the parameter space is too
large, the essential
difficulty of such problems cannot be captured by testing a simple
null versus a simple alternative.
Instead rate optimal lower bounds can often be provided by testing a simple
null versus a composite alternative where
the value of the functional is constant on the composite alternative.
See, for example, \citet{CL05}, where upper and lower bounds are
constructed for quadratic functionals over many different parameter spaces.

Recently some nonsmooth functionals have been considered. A
particularly interesting paper is \citet{LNS99} which studies the problem
of estimating the $L_r$ norm of the drift function under the white
noise model. One of the key observations in this paper is the need to
consider testing between two
composite hypotheses where the $L_r$ norm is not
constant on either of these composite hypotheses and where the sets of
values of the functional on these two hypotheses are interwoven. These
are called fuzzy hypotheses in the language of \citet{Tsybakov09}.

The purpose of the present paper is to advance these ideas further. We
first develop a new general minimax lower bound technique for
estimating any functional $T$ based on testing two composite
hypotheses. For any two priors, say $\mu_0$ and $\mu_1$, on the
parameter space we obtain a lower bound on the expected squared bias
with respect to $\mu_1$ under a constraint on the upper bound of the
expected mean squared error with respect to $\mu_0$. The lower bound
depends on the difference between the expected value of $T$ over each
of the priors and also on the variance of $T$ under~$\mu_0$. The bound
also depends on the chi-square distance between the two marginal
distributions of the observations, one over $\mu_0$, the other over
$\mu_1$. Some of the technical tools for deriving minimax lower bounds
developed earlier in the literature can be seen as special cases of the
general result given in the present paper.

We then consider specifically the problem of estimating the $\ell_1$
norm of a multivariate normal mean vector. This nonsmooth functional
estimation problem exhibits some features that are significantly
different from those in estimating smooth functionals in terms
of the optimal rates of convergence as well as the technical tools
needed for the analysis of both the minimax lower bounds and the
construction of the optimal estimators.

Let $y_1, y_2,\ldots, y_n$ be independent normal random variables
where $y_i \sim N(\theta_i, 1)$.
The problem of focus in this paper is that of estimating
%
%
\begin{equation}
\label{T}
T(\theta) = {1 \over n} \sum_{i=1}^n |\theta_i|,
\end{equation}
where we assume
that either $|\theta_i| \le M$ for some constant $M > 0 $ or
that there are no constraints on the $\theta_i$.
In the present paper we develop optimal estimators of $T(\theta)$
along with minimax lower bounds.
In particular for the bounded case we construct an asymptotically
sharp minimax estimator using approximation theory and
Hermite polynomials.
By combining the minimax lower and upper
bounds developed in later sections, the main results on the minimax
estimation of the functional $T(\theta)$
can be summarized in the following theorem.
\begin{theorem}
Let $Y \sim N(\theta, I_n)$ and let
$T(\theta) = {{1 \over n} \sum_{i=1}^n} |\theta_i|$. For a fixed
constant $M> 0$, denote by
$\Theta_n(M) = \{\theta\in\RR^n\dvtx|\theta_i|\le M\}$.
Then the minimax risk for estimating the functional $T(\theta)$ based
on $Y$ over $\Theta_n(M)$ satisfies
%
%
\begin{equation}
\inf_{\hat T} \sup_{\theta\in\Theta_n(M)} E\bigl(\hat T -
T(\theta)\bigr)^2= \beta_*^2 M^2 \biggl({\log\log n \over\log n}
\biggr)^2\bigl(1+o(1)\bigr),
\end{equation}
where $\beta_*\approx0.28017$ is the Bernstein constant, and the
minimax risk
for estimating the functional $T(\theta)$ over $\RR^n$ satisfies
%
%
\begin{equation}
\inf_{\hat T} \sup_{\theta\in\RR^n} E\bigl(\hat T - T(\theta)\bigr)^2
\asymp{1 \over\log n}.
\end{equation}
\end{theorem}

These rates are dramatically different from the usual parametric or
algebraic rates of convergence for estimating smooth functionals. The
fundamental difficulty of estimating the functional $T(\theta)$ can be
traced back
to the nondifferentiability of the absolute value function at the
origin. This is reflected both in the derivation of the lower bounds
and the construction of the optimal estimators. Best polynomial approximation
and Hermite polynomials play major roles in the derivation of the lower
bounds as well as in the construction of the optimal estimators.

The minimax lower bounds are established by applying the
general lower bound technique to two carefully
constructed composite hypotheses.
In the present context to obtain good lower bounds,
neither prior can be degenerate.
A key step is the construction of two mixture priors which have a
large difference in the expected values of
the functional while making the chi-square distance between the two
mixture models small. In order to turn this heuristic idea into an
effective tool it is necessary to be able to bound the chi-square
distance between two normal mixture models. In previous applications
such bounds have only been given in the much simpler case when one of
the mixtures is degenerate. See, for example, \citet{CL05} and
\citet{WBCL08}.

The construction of the optimal estimators of the nonsmooth
functional $T(\theta)$ is significantly more complicated than those for
linear or quadratic functionals. For optimal estimation of $T(\theta)$
over the bounded set $\Theta_n(M)$, we first use the best polynomial
approximation $G_K^*(x)= \sum_{k=0}^{K} g_{2k}^* x^{2k}$ of the absolute
value function~$|x|$. Then for each $i$ and each $k$ we form an
unbiased estimate of $\theta_i^k$ using the Hermite
polynomials. Putting these terms together for a given $i$ yields an
estimate of $|\theta_i|$. An effective estimate of the functional $T$
can then be constructed by averaging these estimates of~$|\theta_i|$.
We show that by carefully selecting the cutoff $K=K_n$
the resulting estimator is asymptotically sharp minimax.
This estimator is, however, not optimal over the unbounded parameter
space $\RR^n$. An additional testing step is used to construct a
hybrid estimator and it is shown that
the estimator is rate optimal for estimating $T(\theta)$ over $\RR^n$.
In addition, we also consider the estimation of $T(\theta)$ over a
parameter space where the mean $\theta$ is a high-dimensional sparse
vector with a small fraction of nonzero coordinates.

The rest of the paper is organized as follows. In Section \ref{lowerbd.sec}
we derive the general lower bounds for estimating any functional $T$
based on testing two composite hypotheses.
In Section \ref{normal.lowerbd.sec} we bound the chi-square distance
between two normal mixture models and apply the general lower bound
from Section \ref{lowerbd.sec} to derive minimax lower bounds for
estimating the nonsmooth functional $T(\theta)$ given in (\ref{T}).
Section \ref{bounded.upperbnd.sec} constructs an estimator of
$T(\theta)$ using best polynomial approximation and Hermite
polynomials and shows that the estimator is sharp minimax for the
bounded case. Section \ref{unbounded.upperbnd.sec} considers the
unbounded case. A hybrid estimator is constructed and is shown to
attain the optimal rate of convergence. Section
\ref{sparse.upperbnd.sec} treats the sparse case.
Discussions on the connections and differences of
our results with other related work is given in Section
\ref{discussion.sec}. Technical lemmas and some of
the main results are proved in Section \ref{proofs.sec}.

\section{General lower bound}
\label{lowerbd.sec}

In this section, a constrained risk
inequality is developed which immediately yields a general minimax
lower bound
based on testing two composite hypotheses.

Suppose we observe a random variable $X$ which has a distribution
$P_\theta$ where $\theta$ belongs to a given parameter
space $\Theta$. Let $\hat T=\hat T(X)$ be an estimator of a function
$T(\theta)$ based on $X$ and denote the bias of $\hat T$ by
$B(\theta) = E_\theta\hat T - T(\theta)$.
Let $\Theta_0$ and $\Theta_1$ be subsets of the parameter space
$\Theta$ where $\Theta_0 \cup\Theta_1 = \Theta$.
Let $\mu_0$ and $\mu_1$ be two prior distributions supported on
$\Theta_0$ and $\Theta_1$, respectively.

Let $m_i$ and $v^2_i$ be the means and variances of $T(\theta)$ under
the priors $\mu_i$ for $i=0$ and $1$.
More specifically,
\[
m_i = \int T(\theta) \mu_i(d\theta) \quad\mbox{and}\quad
v_i^2 = \int\bigl(T(\theta) - m_i\bigr)^2 \mu_i(d\theta).
\]
Write $F_i$ for the marginal distribution of $X$ when the prior is
$\mu_i$ for $i=0,1$. Let $f_i$ be the density of $F_i$ with respect to
a common dominating measure of $F_0$ and~$F_1$.
For any function $g$ we shall write $E_{f_{0}}g(X)$ for the expectation of
$g(X)$ with respect to the marginal distribution of $X$ when the prior on
$\theta$ is $\mu_o$. We shall write
$E_{\theta}g(X)$ for the expectation of $g(X)$ under $P_{\theta}$.

Finally define
the chi-square distance between $f_0$ and $f_1$ by
\[
I = \biggl\{E_{f_{0}} \biggl({f_1(X) \over f_0(X)} - 1 \biggr)^2
\biggr\}^{1/2}.
\]

The following theorem gives a lower bound for the average risk of an estimator
$\hat T$ under any mixture prior $\lambda\mu_0 + (1-\lambda)\mu_1$,
$0\le\lambda\le1$.
\begin{theorem}
\label{lower.bd.thm}
{\smallskipamount=0pt
\begin{longlist}
\item
Suppose
$\int E_{\theta}(\hat T(X) - T(\theta))^2 \mu_0(d\theta) \le
\varepsilon^2$,
then
%
%
\begin{equation}
\label{lb1}
\biggl| \int B(\theta) \mu_1(d\theta) - \int B(\theta) \mu
_0(d\theta) \biggr|
\ge|m_1 - m_0| - (\varepsilon+ v_0)I .
\end{equation}
\item
If $|m_1 - m_0|> v_0I$
and
$0 \le\lambda\le1$, then
%
%
\begin{eqnarray}
\label{bayes}
&&\int E_{\theta}\bigl(\hat T(X) - T(\theta)\bigr)^2 \bigl(\lambda\mu_0(d\theta) +
(1-\lambda)
\mu_1(d\theta)\bigr) \nonumber\\[-8pt]\\[-8pt]
&&\qquad\ge
{\lambda(1-\lambda)(|m_1-m_0| - v_0I)^2 \over\lambda+ (1
-\lambda)(I+1)^2}\nonumber
\end{eqnarray}
and in particular
%
%
\begin{equation}
\label{minimax}
\max_{i=0,1} \int E_{\theta}\bigl(\hat T(X) - T(\theta)\bigr)^2 \mu
_i(d\theta) \ge
{(|m_1-m_0| - v_0I)^2 \over(I+2)^2}.
\end{equation}
\end{longlist}}
\end{theorem}

Informally, Theorem \ref{lower.bd.thm} says that if the average risk of
$\hat T$ under $\mu_0$ is ``small,'' then the change in
average bias under $\mu_0$ and under $\mu_1$ must be ``large.''
In particular, this implies that the average risk under a mixture
prior is ``large.''

Since the maximum risk is always at least as large as the average
risk, Theorem~\ref{lower.bd.thm} yields immediately a lower bound on
the minimax risk.
\begin{corollary}
\label{minimax.cor}
If $|m_1 - m_0| > v_0 I$, then
%
%
\begin{equation}
\label{minimax2}
\sup_{\theta\in\Theta} E_{\theta}\bigl(\hat T(X) - T(\theta)\bigr)^2 \ge
{(|m_1-m_0| - v_0I)^2 \over(I+2)^2}.
\end{equation}
\end{corollary}

Simpler versions of constrained risk inequalities have been developed
before, most often for studying the cost of adaptation and
superefficiency.
For example, a~two-point risk inequality was given in \citet{BL96}
and used to study adaptive estimation of linear functionals.
The constrained risk inequality given in the present paper allows for a richer
collection of applications and is especially useful
when estimating nonsmooth functionals where it is essential
to test complicated composite hypotheses in order to obtain good minimax
lower bounds.
In particular the lower bounds given in the next section rely on
Corollary \ref{minimax.cor}.
\begin{pf*}{Proof of Theorem \ref{lower.bd.thm}}
We shall also assume without loss of generality that $m_1 \ge m_0$.
Then
\begin{eqnarray*}
&&E_{f_0}\biggl\{\bigl(\hat T(X) - m_0\bigr)\biggl({f_1(X) -f_0(X) \over f_0(X)}\biggr)
\biggr\} \\
&&\qquad=
m_1 + \int
B(\theta) \mu_1(d\theta) - \biggl(m_0 +
\int
B(\theta) \mu_0(d\theta) \biggr).
\end{eqnarray*}
Now note that
\begin{eqnarray*}
&&
E_{f_0}\bigl(\hat T(X) - m_0\bigr)^2 \\
&&\qquad= \int E_{\theta}\bigl(\hat T(X) - m_0\bigr)^2
\mu_0(d\theta)\\
&&\qquad= \int E_{\theta}\bigl(\hat T(X) - T(\theta) + T(\theta) - m_0 \bigr)^2
\mu_0(d\theta)\\
&&\qquad= \int E_{\theta}\bigl(\hat T(X) - T(\theta)\bigr)^2
\mu_0(d\theta)\\
&&\qquad\quad{} + \int\bigl(T(\theta) - m_0 \bigr)^2\mu_0(d\theta)\\
&&\qquad\quad{}+ 2\int B(\theta) \bigl(T(\theta) - m_0\bigr) \mu_0(d\theta) \\
&&\qquad \le \varepsilon^2 + v_0^2 + 2v_0\varepsilon= (\varepsilon+v_0)^2.
\end{eqnarray*}
The Cauchy--Schwarz inequality now yields
\begin{eqnarray*}
E_{f_0}\biggl\{\bigl(\hat T(X) - m_0\bigr)\biggl({f_1(X) -f_0(X) \over f_0(X)}\biggr)\biggr\}
&\le&\bigl(E_{f_0}\bigl(\hat T(X) - m_0\bigr)^2\bigr)^{1/2}\cdot I \\
&\le&(\varepsilon
+ v_0) I.
\end{eqnarray*}
Hence,
%
%
\begin{equation}
m_1 + \int
B(\theta) \mu_1(d\theta) - \biggl(m_0 +
\int
B(\theta) \mu_0(d\theta) \biggr)
\le(\varepsilon+ v_0) I,
\end{equation}
and it follows that
\[
\int B(\theta) \mu_1(d\theta) - \int B(\theta) \mu_0(d\theta)
\le m_0 - m_1 + (\varepsilon+ v_0) I,
\]
which in turn yields (\ref{lb1}).

Now consider the quadratic
%
%
\begin{equation}
\label{eq:2}
J(x) = \lambda x^2 + (1-\lambda)(a -bx)^2,
\end{equation}
where we assume that $0 < \lambda<1$, $a>0$ and $b>0$.
It is easy to check that $J$ is minimized when $x = x_{\min
}={ab(1-\lambda) \over
\lambda+ b^2(1 -\lambda)}$
and that at this value $a-bx >0$
and $J(x_{\min}) = {a^2\lambda(1-\lambda) \over{\lambda+
b^2(1-\lambda)}}$.
It follows that
%
%
\begin{equation}
\lambda x^2 + (1-\lambda)\bigl( \max(a -bx,0)\bigr)^2
\end{equation}
is also minimized at this same value.
Now we also have
\[
\int B^2(\theta) \mu_1(d\theta)
\ge
\bigl(\max\bigl( m_1 - m_0 - v_0 I - (I+1)\varepsilon,0\bigr)\bigr)^2.
\]
It follows that for $0 \le\lambda\le1$
\begin{eqnarray*}
&&\lambda\varepsilon^2 + (1 - \lambda) \int B^2(\theta)
\mu_1(d\theta)\\
&&\qquad \ge \lambda\varepsilon^2 + (1 - \lambda) \bigl(\max\bigl( m_1 - m_0 - v_0 I
- (I+1)\varepsilon,0\bigr)\bigr)^2\\
&&\qquad\ge {\lambda(1-\lambda)(|m_1-m_0| -v_0I)^2 \over{\lambda+
(1-\lambda)(I+1)^2}},
\end{eqnarray*}
which gives (\ref{bayes}). The final inequality (\ref{minimax})
follows by setting $\lambda= {I+1 \over I+2}$
since the minimax risk is greater than any Bayes risk.
\end{pf*}

\section{Lower bound for estimating the $\ell_1$ norm of normal means}
\label{normal.lowerbd.sec}

We now turn to the problem of optimally estimating a particular nonsmooth
functional where the use of the lower bound developed in the previous
section yields sharp results. Let
$y_i \stackrel{\mathrm{ind}}{\sim} N(\theta_i, 1)$, $i=1, 2,\ldots, n$,
and consider the functional $T$ where
%
%
\begin{equation}
T(\theta) = {1\over n} \sum_{i=1}^n |\theta_i|.
\end{equation}
As mentioned in the \hyperref[intro]{Introduction}, there are two particularly
interesting cases. One is the bounded case with $\theta\in\Theta
_n(M)$ where
$\Theta_n(M) = \{\theta\in\RR^n\dvtx\break|\theta_i|\le M\}$ with a
constant $M>0$.
Another case is the unbounded case where $\theta\in\RR^n$.
It is worth noting that we need to consider the bounded case
with a bound growing in $n$ in order to solve the unbounded case.
In addition, we are also interested in the sparse case where $\theta$
is a high-dimensional sparse vector with a small fraction of nonzero
coordinates.

In this section the focus is on developing minimax lower bounds. The
minimax upper bounds and the optimal estimation procedures will be
given in the next three sections. Best polynomial approximation
plays a major role in the development of the lower bound and as we
shall see later also in the development of the upper bound.

\subsection{Best polynomial approximation of the absolute value function}

Optimal polynomial approximation of the absolute value function has
been well studied in approximation theory. See, for example, \citet
{Bernstein13}, \citet{VC87} and \citet{Rivlin90}. For a given
positive integer $k$,
let $\mathcal{P}_k$ denote the class of all real polynomials
of degree at most $k$. For any continuous function $f$ on $[-1, 1]$, let
\[
\delta_k(f) = {\inf_{G\in\mathcal{P}_k} \max_{x\in[-1,1]}} |f(x) - G(x)|.
\]
A polynomial $G^*$ is said to be a best polynomial approximation of
$f$ if
\[
\delta_k(f) = {\max_{x\in[-1,1]} }|f(x) - G^*(x)|.
\]
We now focus on the special case of the absolute value function
$f(x)=|x|$. Because $f$ is an even function, so is its best polynomial
approximation. We thus only need to consider polynomials of even
degrees. For any positive integer $k$, we shall denote by $G_k^*$ the
best polynomial approximation of degree $2k$ to $|x|$ and write
%
%
\begin{equation}
G_k^*(x)=\sum_{j=0}^k g_{2j}^* x^{2j} .
\end{equation}
The Bernstein constant is defined as
\[
\beta_* =\lim_{k \to\infty}2k\delta_{2k}(f).
\]
\citet{Bernstein13} showed that the limit exists and is between
0.278 and
0.286. \citet{VC87} disproved a conjecture by Bernstein
and calculated that $\beta_* = 0.280169499$.

The classical Chebyshev alternation theorem states that a polynomial
$G^*\in\mathcal{P}_k$ is the (unique) best polynomial approximation
to a
continuous function $f$ if and only if
the difference $f(x) - G^*(x)$ takes consecutively its maximal value with
alternating signs at least $(k + 2)$ times. That is,
there exist $k + 2$ points $-1\le x_0 < \cdots< x_{k+1} \le1$ such that
\[
[f(x_j) - G^*(x_j)] = \pm(-1)^j \max_{x\in[-1,1]} |f(x) -
G^*(x)|,\qquad
j= 0, \ldots, k + 1.
\]
In the case of the absolute value function, the best polynomial
approximation $G_k^*(x)$ has at least $2k+2$ alternation
points. The set of these alternation points is important in the
construction of the least favorable priors used in the derivation of
the minimax lower bounds given in this section.
Divide the set of the alternation points of $G_k^*(x)$ into two
subsets and denote
%
%
\begin{eqnarray}
\label{A0}
A_0&=&\{x\in[-1, 1]\dvtx|x| - G_k^*(x) = -\delta_{2k}(|x|)\},\\
\label{A1}
A_1&=&\{x\in[-1, 1]\dvtx|x| - G_k^*(x) = \delta_{2k}(|x|)\}.
\end{eqnarray}
It follows easily from the fact that both $|x|$ and $G_k^*(x)$ are even
functions that the set $A_0$ contains an odd number of points and $A_1$
has an even number of points. We shall see later that least
favorable priors are necessarily supported on $A_0$ and~$A_1$, respectively.
Intuitively, this makes the priors maximally apart and yet not
``testable.'' It also connects the construction of the optimal
estimator with the minimax lower bound.

\subsection{Minimax lower bounds}

We now state and prove the minimax lower bounds for estimating the
nonsmooth functional $T(\theta)$ over the bounded set $\Theta_n(M)$
and the unbounded set $\RR^n$. The derivation of the lower bounds
relies heavily on the general lower bound argument given in the previous
section. It also requires a careful construction of least
favorable prior distributions $\mu_0$ and $\mu_1$ along with finding
an effective upper bound for the chi-square distance between the
marginal distributions.
\begin{theorem}
\label{normal.lower.bd.thm}
Let\vspace*{1pt} $y_i\sim N(\theta_i, 1)$, $i=1,\ldots, n$, be independent normal
random variables, and let $T(\theta) = {{1 \over n} \sum_{i=1}^n}
|\theta_i|$.
For a\vspace*{1pt} fixed constant $M> 0$, denote by
$\Theta_n(M) = \{\theta\in\RR^n\dvtx|\theta_i|\le M\}$.
Then, the minimax risk for estimating $T(\theta)$ over the parameter
space $\Theta_n(M)$ is bounded from below as
%
%
\begin{equation}
\label{fixed}
\inf_{\hat T} \sup_{\theta\in\Theta_n(M)} E\bigl(\hat T - T(\theta)\bigr)^2
\ge\beta_*^2 M^2 \biggl({\log\log n \over\log n}\biggr)^2\bigl(1+ o(1)\bigr),
\end{equation}
where $\beta_*$ is the Bernstein constant.
Without any constraint on the parameters, the minimax risk satisfies
%
%
\begin{equation}
\label{unconstrained}
\inf_{\hat T} \sup_{\theta\in\RR^n} E\bigl(\hat T - T(\theta)\bigr)^2
\ge{4\beta_*^2 \over9e^2\log n}\bigl(1+ o(1)\bigr).
\end{equation}
\end{theorem}

The minimax lower bounds given in Theorem \ref{normal.lower.bd.thm}
converge to zero at a slow logarithmic rate showing
that the nonsmooth functional $T(\theta)$ is difficult to estimate.
In contrast the rates for estimating
linear and quadratic functionals are most often algebraic. In
particular let
\[
L(\theta) = {1 \over n} \sum_{i=1}^n \theta_i \quad\mbox{and}\quad
Q(\theta) = {1 \over n} \sum_{i=1}^n \theta_i^2.
\]
It is easy to check that the usual parametric rate of
convergence over $\RR^n$ for estimating the linear functional
$L(\theta)$ can be attained by the sample average~$\bar y$. For
estimating the quadratic functional
$Q(\theta)$, the parametric rate can be achieved over
$\Theta_n(M)$ by using the unbiased estimator
$\hat Q = {1\over n}\sum_{i=1}^n (y_i^2 - 1)$.

We shall show in the next section that the minimax lower bound\break
$\beta_*^2M^2({\log\log n \over\log n})^2$ for
$\Theta_n(M)$ is in fact asymptotically sharp and
the rate of convergence ${1\over\log n}$ for $\RR^n$ is optimal. The
optimal procedures are constructed using the Hermite
polynomials. These procedures are much more involved than those for
estimating the linear and quadratic functionals discussed above.



A crucial tool in the proof of the lower bounds as well as in the
construction of the optimal procedures is the application of
properties of Hermite polynomials.
Let $H_k$ be the Hermite polynomial defined by
%
%
\begin{equation}
\label{hermite}
{d^k \over dy^k} \phi(y) = (-1)^kH_k(y)\phi(y).
\end{equation}
For this version of the Hermite polynomial
%
%
\begin{equation}
\int H_k^2(y) \phi(y) \,dy = k! \quad\mbox{and}\quad
\int H_k(y)H_j(y) \phi(y) \,dy = 0,
\end{equation}
when $k \neq j$.

Another key technical tool for the proof of Theorem \ref{normal.lower.bd.thm}
is the construction of two priors with special properties.
\begin{lemma}
\label{prior.lem}
For any given even integer $k>0$, there exist two probability
measures $\nu_0$ and $\nu_1$ on $[-1,1]$ that satisfy the following
conditions:
\begin{itemize}
\item$\nu_0$ and $\nu_1$ are symmetric around $0$;
\item
$\int t^l \nu_1(dt) =\int t^l \nu_0(dt) $, for $l=0,1, \ldots, k$;
\item$\int|t| \nu_1(dt) - \int|t| \nu_0(dt)
= 2 \delta_{k}$,
\end{itemize}
where $\delta_{k}$ is the distance in the uniform norm on $[-1, 1]$
from the absolute value function $f(x)=|x|$ to the space
$\mathcal{P}_k$ of
polynomials of no more than degree $k$.
\end{lemma}

As discussed earlier, $\delta_{k} = \beta_* k^{-1}(1+o(1))$
as $k\goto\infty$, where $\beta_*$ is the Bernstein constant.
See Section \ref{discussion.sec} for further discussions. The proof of
Lemma \ref{prior.lem} is given in Section \ref{proofs.sec}.
\begin{pf*}{Proof of Theorem \ref{normal.lower.bd.thm}}
For a given even integer $k_n$, let $\nu_0$ and $\nu_1$ be two
probability measures possessing the properties given in Lemma \ref{prior.lem}.
Let $g(x) = M x$ and let $\mu_i$ be the measures on $[-M, M]$
defined by $\mu_i (A) = \nu_i (g^{-1}(A))$ for $i=0$ and 1.
It follows that:
\begin{itemize}
\item$\mu_0$ and $\mu_1$ are symmetric around $0$;
\item$\int t^l \mu_1(dt) =\int t^l \mu_0(dt) $, for $l=0,1, \ldots
, k_n$;
\item$\int|t| \mu_1(dt) - \int|t| \mu_0(dt)
= 2M\delta_{k_n}$.
\end{itemize}

Let $\mu_1^n$ and $\mu_0^n$ be the product priors
$\mu_i^n = \prod_{j=1}^n \mu_i$.
In other words, we put down $n$ independent priors on the coordinates.
We have
\[
E_{\mu_1^n} T(\theta)
- E_{\mu_0^n} T(\theta)
= E_{\mu_1} |\theta_1| -E_{\mu_0} |\theta_1|
= 2M\delta_{k_n}
\]
and
\[
E_{\mu_0^n} \bigl(T(\theta) - E_{\mu_0^n} T(\theta) \bigr)^2
= {1 \over n} E_{\mu_0} (|\theta_1| - E_{\mu_0} |\theta_1| )^2
\le{ M^2 \over n}.
\]
Set $f_{0,M}(y) = \int\phi(y- t)\mu_0(dt)$
and $f_{1,M}(y) = \int\phi(y-t)\mu_1(dt)$.
Note that since $g(x) = \exp(-x)$ is a convex function of $x$,
and $\mu_0$ is symmetric,
\begin{eqnarray*}
f_{0,M}(y) & \ge& {1 \over{\sqrt{2 \pi}}}\exp\biggl(-\int{(y - t)^2
\over2}\mu_0(dt)\biggr)\\
& = &\phi(y) \exp\biggl(-{1\over2}M^2\int t^2 \nu_0(dt)\biggr)\\
& \ge& \phi(y) \exp\biggl(-{1\over2}M^2\biggr).
\end{eqnarray*}
Let $H_r$ be the Hermite polynomial defined in (\ref{hermite}).
Then
\[
\phi(y-\alpha t) = \sum_{k=0}^{\infty} H_k(y) \phi(y) {\alpha^k t^k
\over{k!}},
\]
and it follows that
\[
\int{(f_{1,M}(y) - f_{0,M}(y))^2 \over f_{0,M}(y) }\,dy
\le
e^{{M^2 /2}}\sum_{k=k_n +1}^{\infty}
{1 \over k!} M^{2k}.
\]
Now set
\[
I^2_n =
\int{(\prod_{i=1}^nf_{1,M}(y_i) -
\prod_{i=1}^nf_{0,M}(y_i))^2 \over\prod_{i=1}^nf_{0,M}(y_i)
}\,dy_1\,dy_2\cdots dy_n.
\]
Then
%
%
\begin{eqnarray}\label{bb}
I_n^2
& =&
\int{(\prod_{i=1}^nf_{1,M}(y_i)
)^2 \over\prod_{i=1}^nf_{0,M}(y_i)
}\,dy_1\,dy_2\cdots dy_n -1 \nonumber\\
& = &
\Biggl(\prod_{i=1}^n\int{(
f_{1,M}(y_i)
)^2 \over f_{0,M}(y_i)
}\,dy_i\Biggr) -1 \nonumber\\[-8pt]\\[-8pt]
& \le& \Biggl(1 + e^{{M^2 /2}}\sum_{k=k_n +1}^{\infty}
{1 \over k!} M^{2k} \Biggr)^n -1
\nonumber\\
& \le& \biggl(1 + e^{{3M^2 /2}}
{1 \over k_n!} M^{2k_n} \biggr)^n -1. \nonumber
\end{eqnarray}
Now note that $k! > ({k \over e})^k$.
Hence
%
%
\begin{equation}
\label{bound}
I_n^2
\le
\biggl(1 + e^{{3M^2 /2}}
\biggl({e M^{2}\over k_n}\biggr)^{k_n} \biggr)^n -1 .
\end{equation}
Now let $k_n$ be the smallest positive integer satisfying
$k_n \ge{\log n \over\log\log n} + {\log n \over(\log\log n)^{3/2}}$.
It is easy to check that $I_n \rightarrow0$.
Noting that $v_0 \le{ M \over\sqrt n}$
and applying Corollary~\ref{minimax.cor} yields
\begin{eqnarray*}
\inf_{\hat T} \sup_{\theta\in\Theta_n(M)} E\bigl(\hat T - T(\theta)\bigr)^2
&\ge&{( 2M\delta_{k_n} - ({M /\sqrt n})I_n)^2
\over(I_n+2)^2}\\
&=& \beta_*^2 M^2 \biggl({\log\log n \over\log n}\biggr)^2 \bigl(1+ o(1)\bigr),
\end{eqnarray*}
and (\ref{fixed}) follows.

For the proof of (\ref{unconstrained}), let $M= \sqrt{\log n}$
and take $k_n$ to be the smallest positive integer satisfying $k_n \ge
(1.5)e \log n$.
We may bound $I_n$ starting from (\ref{bb}) and then noting that for
some constant $D>0$
%
%
\begin{eqnarray}
I_n^2
&\le&
\biggl(1 + e^{{M^2 /2}}D
{1 \over k_n!} M^{2k_n} \biggr)^n -1\nonumber\\[-8pt]\\[-8pt]
&\le&\biggl(1 + Dn^{1 /2} \biggl({e \log n \over{(3/2)}e \log n}\biggr)^{k_n}\biggr)^n -1
\rightarrow0.\nonumber
\end{eqnarray}
It is then easy to check that Corollary \ref{minimax.cor} now yields
(\ref{unconstrained}).
\end{pf*}
\begin{remark}
In the bounded case, we shall show in Section \ref{bounded.upperbnd.sec}
that the minimax lower bound
$\beta_*^2M^2({\log\log n \over\log n})^2$ is
asymptotically sharp. It can be seen from the proof of
Theorem \ref{lower.bd.thm} that this minimax risk corresponds to
the Bayes risk of the least favorable prior which is asymptotically
equal to the prior ${1 \over2} (\mu_0 + \mu_1)$.
\end{remark}
\begin{remark}
The proof of (\ref{unconstrained}) can be used to show that
for any constant $c>0$, there exists another constant $d>0$ such that
%
%
\begin{equation}
\label{growing.lowerbd}
\inf_{\hat T} \sup_{\theta\in\Theta_n(\sqrt{c\log n})}
E\bigl(\hat T - T(\theta)\bigr)^2 \ge{d \over\log n}\bigl(1+ o(1)\bigr).
\end{equation}
\end{remark}

\section{Optimal estimation of the $\ell_1$ norm of bounded normal means}
\label{bounded.upperbnd.sec}

Section~\ref{normal.lowerbd.sec} developed minimax lower bounds for
estimating the nonsmooth functional $T(\theta)$.
Although the minimax lower bounds converge slowly, they are also
difficult to attain. The difficulty of the estimation problem
stems from the fact that the absolute value function is not
differentiable at $0$. In this section we shall consider
the bounded case and construct an estimator that relies on the best
polynomial approximation to the absolute value function and the use of
Hermite polynomials. The estimator is then shown to be asymptotically
sharp minimax. The unbounded case and the sparse case will be treated
in the next two sections.

\subsection{Polynomial approximation}

The construction of the rate optimal estimator is involved.
This is partly due to the nonexistence of an unbiased estimator
for~$|\theta_i|$.
Our strategy is to ``smooth'' the singularity at $0$ by a~polynomial
approximation and construct an unbiased estimator for
each term in the expansion by using the Hermite polynomials.

The optimal estimator relies on the best polynomial approximation
$G_K^*$ of the absolute value function. A drawback of using $G_K^*$
is that it is not convenient to construct. An explicit and nearly
optimal polynomial approximation $G_K$ can be easily obtained by using
the Chebyshev polynomials. Note that the Chebyshev polynomial (of the
first kind) of degree $k$ is defined
as
\[
T_{k}(x) = \sum_{j=0}^{[{k/2}]} (-1)^j {k\over k - j} \pmatrix{k-j\cr
j} 2^{k-2j-1}
x^{k-2j}.
\]
The following expansion can be found, for example, in Rivlin (\citeyear{Rivlin90}):
%
%
\begin{equation}
\label{infinite.expansion}
|x| = {2 \over\pi}T_0(x) +{4 \over\pi}\sum_{k=1}^{\infty} (-1)^{k+1}
{T_{2k}(x) \over{4k^2-1}},
\end{equation}
where $T_{2k}(x)$ is the Chebyshev polynomial of degree $2k$.
Consider the truncated version of the expansion
(\ref{infinite.expansion}) and let
%
%
\begin{equation}
\label{GK}
G_K(x) = {2 \over\pi}T_0(x) +{4 \over\pi}\sum_{k=1}^{K} (-1)^{k+1}
{T_{2k}(x) \over{4k^2-1}}.
\end{equation}
We can also write $G_K(x)$ as
%
%
\begin{equation}
G_K(x)=\sum_{k=0}^K g_{2k} x^{2k}.
\end{equation}

The following lemma provides uniform error bounds of $G_K^*$ and
$G_K$ over the interval $[-1, 1]$ as well as bounds on the
coefficients $g_{2k}^*$ and $g_{2k}$. These bounds are useful in the
analysis of the optimal estimators.
\begin{lemma}
\label{poly.lem}
Let $G_K^*(x)=\sum_{k=0}^K g_{2k}^* x^{2k}$ be the best polynomial
approximation of degree $2K$ to $|x|$, and let $G_K$ be defined in
(\ref{GK}). Then
%
%
\begin{eqnarray}
\label{Bernstein.bnd*}
\max_{x\in[-1,1]} \bigl| G_K^*(x) - |x| \bigr| &\le&
{\beta_*\over2K}\bigl(1+o(1)\bigr),
\\
\label{Bernstein.bnd}
\max_{x\in[-1,1]} \bigl| G_K(x) - |x| \bigr| & \le& {2 \over\pi(2K+1)}.
\end{eqnarray}
The coefficients $g_{2k}^*$ and $g_{2k}$ satisfy for all $0\le k \le K$,
%
%
\begin{equation}
|g_{2k}^*| \le2^{3K} \quad\mbox{and}\quad |g_{2k}| \le2^{3K}.
\end{equation}
\end{lemma}

The uniform error bounds (\ref{Bernstein.bnd*}) and (\ref{Bernstein.bnd})
were proved in \citet{Bernstein13}. The proof of the bound on the
coefficients $g_{2k}^*$ and $g_{2k}$ is given in Section \ref{proofs.sec}.

\subsection{The construction of the optimal estimator}

We shall now use the best polynomial approximation $G_K^*(x)$ and the
Hermite polynomials to construct an estimator of $T(\theta)$ that is
asymptotically sharp minimax over the bounded parameter space
$\Theta_n(M)$.
We first consider the special case of $M=1$. The case of a general $M$
involves an additional rescaling step.

When $M=1$, it follows from Lemma \ref{poly.lem} that each
$|\theta_i|$ can be well approximated by
$G_K^*(\theta_i) = \sum_{k=0}^K g_{2k}^* \theta_i^{2k}$ on the interval
$[-1, 1]$ and hence the functional
$T(\theta)={1\over n} \sum_{i=1}^n |\theta_i|$ can be approximated by
\[
\tilde T(\theta) = {1\over n} \sum_{i=1}^n G_K^*(\theta_i) =
\sum_{k=0}^K g_{2k}^* b_{2k}(\theta),
\]
where $b_{2k}(\theta)\equiv{1\over n} \sum_{i=1}^n \theta_i^{2k}$.
Note that $\tilde T(\theta)$ is a smooth functional, and we shall
estimate $b_{2k}(\theta)$ separately for each $k$ by using the Hermite
polynomials.

Let $\phi$ be the density function of a standard normal variable.
Recall that for positive integers $k$,
%
%
\begin{equation}
{d^k \over dy^k}\phi(y) = (-1)^k H_k(y) \phi(y),
\end{equation}
where $H_k$ is a Hermite polynomial with respect to $\phi$.
It is well known that for $X\sim N(\mu, 1)$, $H_k(X)$ is an unbiased
estimate of $\mu^k$ for any positive integer~$k$, that is,
$EH_k(X)=\mu^k$.

Since $H_k(y_i)$ is an unbiased estimate of $\theta_i^k$ for each $i$,
we can estimate $b_k(\theta)\equiv{1\over n} \sum_{i=1}^n
\theta_i^k$ by $\bar B_k = {1\over n} \sum_{i=1}^n H_k(y_i)$
and define the estimator of $T(\theta)$~by
%
%
\begin{equation}
\label{seq.est}
\widehat{T_K(\theta)} = \sum_{k=0}^K g_{2k}^* \bar B_{2k}.
\end{equation}

For estimating the functional $T(\theta)$ over the bounded
parameter space $\Theta_n(M)$ for a general $M>0$, we shall first
rescale each $\theta_i$ and then approximate $|\theta_i|$ term by term.
More specifically, let $|\theta'_i| = M^{-1}\theta_i$. Then
$|\theta'_i| \le1$ for $i=1,\ldots,n$ and
\[
\bigl||\theta_i'| - G_K^*(\theta_i')\bigr|\le{\beta_*\over2K}
\bigl(1+o(1)\bigr)\qquad
\mbox{for all $|\theta'_i|\le1$.}
\]
Hence,
\[
\bigl||\theta_i| - \tilde G_K^*(\theta_i)\bigr|\le{\beta_*M\over2K}
\bigl(1+o(1)\bigr)\qquad
\mbox{for all $|\theta_i|\le M$},
\]
where $\tilde G_K^*(x)=\sum_{k=0}^K \tilde g_{2k}^* x^{2k}$ with
$\tilde g_{2k}^* = g_{2k}^* M^{-2k+1}$.

Again, $H_k(y_i)$ is an unbiased estimate of $\theta_i^k$.
We estimate $b_{2k}(\theta)\equiv\break{1\over n} \sum_{i=1}^n \theta
_i^{2k}$ by
%
%
\begin{equation}
\bar B_{2k} = {1\over n} \sum_{i=1}^n H_{2k}(y_i)
\end{equation}
and define the estimator of $T(\theta)$ by
%
%
\begin{equation}
\label{seq.est1}
\widehat{T_{K}(\theta; M)} = \sum_{k=0}^{K} \tilde g_{2k}^* \bar B_{2k}
=\sum_{k=0}^{K} g_{2k}^* M^{-2k+1} \bar B_{2k}.
\end{equation}
The performance of the estimator $\widehat{T_{K}(\theta; M)}$ clearly
depends on the choice of the cutoff $K$. We shall specifically choose
%
%
\begin{equation}
\label{K*}
K= K_* \equiv{\log n \over2 \log\log n}
\end{equation}
and define our final estimator of $T(\theta)$ by
%
%
\begin{equation}
\label{bounded.seq.est}
\widehat{T_*(\theta)} \equiv\widehat{T_{K_*}(\theta; M)} = \sum
_{k=0}^{K_*} \tilde g_{2k}^* \bar B_{2k}.
\end{equation}

\subsection{Optimality of the estimator}

We now study the property of the estimator defined in (\ref{bounded.seq.est}).
The following result shows that the estimator
$\widehat{T_{*}(\theta)}$ is asymptotically sharp minimax, that is, it
achieves the exact minimax lower bound given in Theorem
\ref{normal.lower.bd.thm} asymptotically.
\begin{theorem}
\label{seq.L1.thm}
Let $y_i\sim N(\theta_i, 1)$ be independent normal random variables
with $|\theta_i|\le M$, $i=1,\ldots,n$. Let $T(\theta) = {n^{-1}
\sum_{i=1}^n }|\theta_i|$.
The estimator $\widehat{T_{*}(\theta)}$ given in (\ref{bounded.seq.est})
satisfies
%
%
\begin{equation}
\sup_{\theta\in\Theta_n(M)}E\bigl(\widehat{T_{*}(\theta)} - T(\theta)\bigr)^2
\le\beta_*^2 M^2 \biggl({\log\log n \over\log n}\biggr)^2 \bigl(1+o(1)\bigr).
\end{equation}
\end{theorem}
\begin{remark}
If $G_K(x)$, instead of $G_K^*(x)$, is used in the
construction of the estimator $\widehat{T_{*}(\theta)}$, the resulting
estimator $\widehat{T(\theta)}$ satisfies
%
%
\begin{equation}
\sup_{\theta\in\Theta_n(M)}E\bigl(\widehat{T(\theta)} - T(\theta)\bigr)^2
\le4\pi^{-2}M^2 \biggl({\log\log n \over\log n}\biggr)^2 \bigl(1+o(1)\bigr).
\end{equation}
The ratio of this upper bound to the minimax risk is
$4\pi^{-2}/\beta_*^2 \approx5.16$.
\end{remark}

We need the following variance bounds for the proof of
Theorem \ref{seq.L1.thm} as well as other results given in the later sections.
\begin{lemma}
\label{mean-var.lem1}
Let $X\sim N(\mu, 1)$, then
\[
E(H_k^2(X)) = k! \sum_{j=0}^k \pmatrix{k\cr j} \mu^{2j} {1\over j!}.
\]
Consequently
\[
\operatorname{Var}(H_k(X)) \le E(H_k^2(X))\le e^{\mu^2} k^k.
\]
If $|\mu|\le M$ and $M^2 \ge k$, then
\[
\operatorname{Var}(H_k(X)) \le E(H_k^2(X))\le(2 M^2)^k.
\]
\end{lemma}

The proof of Lemma \ref{mean-var.lem1} is given in Section \ref{proofs.sec}.
\begin{pf*}{Proof of Theorem \ref{seq.L1.thm}} In the proof we shall assume
$M\ge1$. The case of $M<1$ is similar.
Note that $E \bar B_{2k} = b_{2k}(\theta)$ for $k \ge0$ and hence,
\[
E \widehat{T_K(\theta; M)} = \sum_{k=0}^K \tilde g_{2k}^*
b_{2k}(\theta)
={1\over n}\sum_{i=1}^n \tilde G_K^*(\theta_i).
\]
The bias of $\widehat{T(\theta)}$ can then be bounded easily as
follows. For any $\theta\in\Theta_n(M)$,
\begin{eqnarray*}
|E \widehat{T_K(\theta; M)} - T(\theta)| &=&
\Biggl\vert{1\over n}\sum_{i=1}^n \tilde
G_K^*(\theta_i)-{1\over n}\sum_{i=1}^n |\theta_i|\Biggr\vert\le
{1\over n} \sum_{i=1}^n \bigl\vert\tilde G_K^*(\theta_i)- |\theta
_i|\bigr\vert\\
&\le&{\beta_* M\over2K} \bigl(1+o(1)\bigr).
\end{eqnarray*}
Now we consider the variance of $\widehat{T_K(\theta; M)}$.
It follows from Lemma \ref{mean-var.lem1} that the variance of $\bar
B_{2k}$ satisfies
\[
\operatorname{Var}(\bar B_{2k}) = n^{-2} \sum_{i=1}^n \operatorname
{Var}(H_{2k}(y_i))
\le e^{M^2} (2k)^{2k} n^{-1}.
\]
To bound the variance of $\widehat{T_K(\theta; M)}$, first note that for
any random variables $X_i$, $i=1,\ldots, n$,
%
%
\begin{equation}\label{inq}
E\Biggl(\sum_{i=1}^n X_i\Biggr)^2 \leq\Biggl(\sum_{i=1}^n (E X_i^2)^{1/2}\Biggr)^2.
\end{equation}
It then follows that for all $\theta\in\Theta_n(M)$,
\begin{eqnarray*}
\operatorname{Var}(\widehat{T_K(\theta; M)}) &\le&
\Biggl\{\sum_{k=1}^K |\tilde g_{2k}^*| \operatorname{Var}^{1/2}(\bar B_{2k})
\Biggr\}^2\\
&\le&\Biggl\{\sum_{k=1}^K |g_{2k}^*| M^{-2k+1} e^{M^2/2} (2k)^{k}
\Biggr\}^2\cdot n^{-1}\\
&\le& 2e^{M^2} 2^{8K} K^{2K} n^{-1}.
\end{eqnarray*}
Hence, the mean squared error of $\widehat{T_K(\theta; M)}$ is
bounded by
%
%
\begin{equation}
\label{MSE1}
E\bigl(\widehat{T_K(\theta; M)} - T(\theta)\bigr)^2 \le
{\beta_*^2 M^2 \over(2K)^2}\bigl(1+o(1)\bigr) + 2e^{M^2} 2^{8K} K^{2K} n^{-1}.
\end{equation}
Now set
\[
K_* = {\log n \over2 \log\log n}.
\]
Then the second term in (\ref{MSE1}) is negligible relative to the
first term and we have, for all $\theta\in\Theta_n(M)$,
\[
E\bigl(\widehat{T_{*}(\theta)} - T(\theta)\bigr)^2 \le\beta_*^2
M^2\biggl({\log\log n \over\log n}\biggr)^2 \bigl(1+o(1)\bigr).
\]
\upqed
\end{pf*}

\section{Estimating the $\ell_1$ norm of unbounded normal means}
\label{unbounded.upperbnd.sec}

We now turn to the unbounded case where no restriction is imposed on
the values of the means $\theta_i$. This case is more difficult than
the bounded case. We shall construct an estimator of
$T(\theta)$ that attains the optimal rate of convergence, but not the
optimal constant, for the unbounded case. In the construction below,
both $G_K^*$ and $G_K$ work. For concreteness, hereafter we shall
focus on using $G_K$ instead of the best polynomial approximation
$G_K^*$.

It turns out that a key step toward solving this
general problem is to understand the estimation problem where the
means are bounded with the bound growing with the sample size $n$. We
shall thus first treat this case and then consider rate-optimal
estimation for the general case.

\subsection{Estimating the $\ell_1$ norm with a growing bound}
\label{growing.sec}

Suppose $y_i \stackrel{\mathrm{ind}}{\sim} N(\theta_i, 1)$, $i=1, 2,\ldots, n$,
where $|\theta_i| \le M_n$ for $i=1,\ldots,n$,
with $M_n = \sqrt{c\log n}$ for some $c> 1$.
As in the last section, we estimate $T(\theta)$ by first rescaling
and define the estimator of $T(\theta)$ by
%
%
\begin{equation}
\label{seq.est2}
\widehat{T_K(\theta; M_n)} = \sum_{k=0}^K \tilde g_{2k} \bar B_{2k},
\end{equation}
where $\tilde g_{2k} = g_{2k} M_n^{-2k+1}$ and
$\bar B_{2k} = {1\over n} \sum_{i=1}^n H_{2k}(y_i)$.

\begin{theorem}
\label{seq.L1.thm2}
Let $y_i\sim N(\theta_i, 1)$ be independent normal random variables
with $|\theta_i|\le M_n$, $i=1,\ldots,n$, where $M_n = \sqrt{c\log n}$
for some $c>1$. Let $T(\theta) = n^{-1}\sum_{i=1}^n |\theta_i|$.
The estimator\vspace*{1pt} $\widehat{T_K(\theta; M_n)}$ given in (\ref{seq.est2}) with
$K={1\over7}\log n -(\log n)^{1/2}$ satisfies
%
%
\begin{equation}
\sup_{\theta\in\Theta_n(M_n)} E\bigl(\widehat{T_K(\theta; M_n)} -
T(\theta)\bigr)^2
\le{49 c\over\pi^2} (\log n)^{-1}\bigl(1+o(1)\bigr).
\end{equation}
\end{theorem}

This upper bound together with the minimax lower bound
(\ref{growing.lowerbd}) show that the estimator $\widehat{T_K(\theta
; M_n)}$
defined in (\ref{seq.est2}) with $K={1\over7}\log n -(\log n)^{1/2}$
is minimax rate optimal in this case.
We shall show that the difficulty of estimating $T(\theta)$ over
$\RR^n$ is essentially the same as estimating over $\Theta_n(M_n)$
with an appropriate choice of $M_n$ of order $\sqrt{\log n}$.
However, the construction of the rate-optimal estimator of $T(\theta)$
over $\RR^n$ is much more complicated.

The proof of Theorem \ref{seq.L1.thm2} is given in Section \ref{proofs.sec}.

\subsection{Rate optimal estimator for the unbounded case}
\label{unbounded.est.sec}

We now turn to the unbounded case.
It is helpful to provide some intuition and motivation before we
formally describe the estimation procedure. Consider the
one-dimensional case. Suppose we observe $X\sim N(\mu, 1)$ and wish to
estimate $|\mu|$. Set $M_n = 8\sqrt{\log n}$.
Let $\mu' = M_n^{-1}\mu$. Then $|\mu'| \le1$ and
\[
\bigl||\mu'| - G_K(\mu')\bigr|\le{2\over\pi(2K+1)}\qquad \mbox{for all
$|\mu'|\le1$.}
\]
Hence,
\[
\bigl||\mu| - \tilde G_K(\mu)\bigr|\le{2 M_n \over\pi(2K+1)}\qquad \mbox
{for all
$|\mu|\le M_n$},
\]
where $\tilde G_K(\mu)=\sum_{k=0}^K \tilde g_{2k} \mu^{2k}$ with
$\tilde g_{2k} =
g_{2k} M_n^{-2k+1}$. Again, $H_k(X)$ is an unbiased estimate of $\mu^k$.
Set $K = {1\over12} \log n$ and
%
%
\begin{equation}
\label{SK}
S_K(x) = \sum_{k=0}^K g_{2k} M_n^{-2k+1} H_{2k}(x).
\end{equation}
We define an estimator of $|\mu|$ by a truncated version of $S_K(X)$,
%
%
\begin{equation}
\label{delta}
\delta(X) = \min\{S_K(X), n\}.
\end{equation}

It is easy to see that
$\delta(X)$ is a good estimate of $|\mu|$ when $|\mu|$ is small. On
the other hand, when $|\mu|$ is large, $\delta(X)$ is no longer a
good estimator of $|\mu|$ because the variance of $\delta(X)$ is very
large. When $|\mu|$ is large, a good estimate of $|\mu|$ is simply $|X|$.
Therefore, for the unbounded case, a good strategy is to estimate
$|\mu|$ by $\delta(X)$ when $|X|$ is not too large and estimate
$|\mu|$ by $|X|$ when $|X|$ is large.

We now formally state the procedure for estimating $T(\theta)$ as follows.
We shall first use the idea of sample splitting. Note that observing
$y_i\sim N(\theta_i, 1)$ is equivalent to observing
$y_{il}\stackrel{\mathrm{i.i.d.}}{\sim} N(\theta_i, 2)$, for $l=1, 2$.
[One can generate $y_{i1}$ and $y_{i2}$ from~$y_i$. Let
$z_i\sim N(0, 1)$ be independent of $y_i$ and set $y_{i1} = y_i + z_i$
and $y_{i2} = y_i - z_i$. Then $y_{il}\stackrel{\mathrm{i.i.d.}}{\sim} N(\theta
_i, 2)$.] Write
$x_{il} = {1\over\sqrt{2}} y_{il}$ for $l=1, 2$ and $i=1,\ldots, n$.
Then $x_{il}\stackrel{\mathrm{i.i.d.}}{\sim} N(\theta_i', 1)$, for $l=1, 2$, with
$\theta_i'=\theta_i/\sqrt{2}$. Estimating $T(\theta)$ based on
$\{y_i\}$ is thus equivalent to estimating $\sqrt{2}T(\theta')$ based
on $\{x_{il}\}$. We shall construct an estimate $\widehat{T(\theta')}$
for $T(\theta')$ and estimate $T(\theta)$ by $\sqrt{2}\widehat
{T(\theta')}$.

We define the estimate of
$T(\theta') = {n^{-1} \sum_{i=1}^n }|\theta_i'|$ by
%
%
\begin{equation}
\label{seq.est3}\qquad
\widehat{T(\theta')} = {1\over n} \sum_{i=1}^n
\bigl\{\delta(x_{i1}) I\bigl(|x_{i2}|\le2\sqrt{2\log n}\bigr)
+ |x_{i1}| I\bigl(|x_{i2}| > 2\sqrt{2\log n}\bigr)\bigr\},
\end{equation}
where $\delta(\cdot)$ is defined in (\ref{SK}) and (\ref{delta}), and
define the estimator of $T(\theta) = n^{-1} \times{\sum_{i=1}^n} |\theta_i|$ by
%
%
\begin{equation}
\label{seq.est4}
\widehat{T(\theta)} = \sqrt{2}\widehat{T(\theta')}.
\end{equation}
Here $|x_{i2}|$ is used to test to size of $\theta_i'$, and based on
the test we use either $\delta(x_{i1})$ or $|x_{i1}|$ to estimate
$|\theta_i'|$.

The following theorem shows that $\widehat{T(\theta)}$ attains the
rate of convergence $(\log n)^{-1}$ over the whole parameter space
$\RR^n$.
\begin{theorem}
\label{unbounded.thm}
The estimator $\widehat{T(\theta)}$ defined in (\ref{seq.est3}) and
(\ref{seq.est4}) satisfies, for all $\theta\in\RR^n$,
%
%
\begin{equation}
E\bigl(\widehat{T(\theta)} - T(\theta)\bigr)^2 \le{C\over\log n} \bigl(1+o(1)\bigr)
\end{equation}
for some constant $C>0$.
\end{theorem}

Together with the minimax lower bound given in Theorem
\ref{normal.lower.bd.thm}, Theorem \ref{unbounded.thm} shows that
the hybrid estimator is rate optimal over the parameter space $\RR^n$.
The proof of Theorem \ref{unbounded.thm} is involved and is given in
Section \ref{proofs.sec}. The key is to analyze the bias and variance
of a single component.

\section{Estimating the $\ell_1$ norm of sparse normal means}
\label{sparse.upperbnd.sec}

In high-dimensional problems, an especially interesting case is when the
mean vector is sparse, that is, only a small proportion of the
$\theta_i$'s are nonzero.
Suppose we observe
$y_i \stackrel{\mathrm{ind}}{\sim} N(\theta_i, 1)$, $i=1, 2,\ldots, n$, where
the mean vector $\theta$ is sparse: only a~small
fraction of components are nonzero, and the locations of the nonzero
components are unknown.

Denote the $\ell_0$ quasi-norm by
$\|\theta\|_0 = \operatorname{Card}(\{i\dvtx \theta_i \neq0\})$.
Fix $k_n$, the collection of vectors with exactly $k_n$ nonzero
entries is
\[
\Theta_{k_n} = \ell_0(k_n) = \{\theta\in\RR^n\dvtx \|\theta\|_0 =
k_n\}.
\]

In this section we consider the problem of estimating the average of
the absolute value of the nonzero means. For $\theta\in\Theta_{k_n}$,
%
%
\begin{equation}
T(\theta) = \operatorname{average}\{|\theta_i|\dvtx \theta_i\neq0\} =
{1\over
k_n} \sum_{i=1}^n |\theta_i|.
\end{equation}

We calibrate the sparsity parameter $k_n$ by $k_n = n^\beta$ for
$0<\beta\le1$. The following result shows that for $0<\beta\le\hf$,
it is not possible to estimate the functional $T(\theta)$ consistently.
\begin{theorem}
\label{sparse.lowerbd.thm}
Let $k_n = n^\beta$. Then for all $0<\beta\le\hf$, the minimax risk
satisfies
%
%
\begin{equation}
\inf_{\widehat{T(\theta)}} \sup_{\theta\in\Theta_{k_n}}
E\bigl(\widehat{T(\theta)} - T(\theta)\bigr)^2 \ge C
\end{equation}
for some constant $C>0$.
\end{theorem}

The proof of Theorem \ref{sparse.lowerbd.thm} is analogous to that of
Theorem 7 in \citet{CL04}, and we omit it here for reasons of
space.

We now turn to the more interesting case where $k_n = n^\beta$ with
$\hf<\beta\le1$. The following result show that the minimax rate of
convergence in this case is $(\log n)^{-1}$.
\begin{theorem}
\label{sparse.thm}
Let $k_n = n^\beta$ for some $\hf< \beta< 1$. Then the minimax risk
for estimating the functional $T(\theta)$ over $\Theta_{k_n}$ satisfies
%
%
\begin{equation}
\inf_{\widehat{T(\theta)}}\sup_{\theta\in\Theta_{k_n}}
E\bigl(\widehat{T(\theta)} - T(\theta)\bigr)^2 \asymp{C\over\log n}.
\end{equation}
\end{theorem}

The proof of the lower bound in Theorem \ref{sparse.thm} is similar to
that of Theorem \ref{normal.lower.bd.thm}. The upper
bound can be attained by a modified version of the estimator
$\widehat{T(\theta)}$ defined in (\ref{seq.est3}) and
(\ref{seq.est4}). The key in the construction is to have estimates
of the individual coordinates that perform well when the coordinates
are zero. This can be achieved by using the polynomial approximation
$G_K(x)$ [or $G_K^*(x)$] without the constant term.

As in Section \ref{unbounded.est.sec} set $K = {1\over12}\log n$, and
define
%
%
\begin{equation}
\label{sparse.SK}
\tilde S_K(x) = \sum_{k=1}^K g_{2k} M_n^{-2k+1} H_{2k}(x).
\end{equation}
Note that here the constant term $g_0$ is excluded.
We then define an estimator of $|\mu|$ by truncating $\tilde S_K(X)$,
%
%
\begin{equation}
\label{sparse.delta}
\tilde\delta(X) = \min\{\tilde S_K(X), n^2\}.
\end{equation}
Note that the bias of the estimator $\tilde\delta(X)$ is much smaller
than the bias of $\delta(X)$ when the mean of $X$ is zero.
As in Section \ref{unbounded.est.sec} we split the sample into two
parts and use one for testing and the other for estimation.
Let $x_{il}$ be defined as in Section~\ref{unbounded.est.sec}. That is,
$x_{il} \stackrel{\mathrm{i.i.d.}}{\sim} N(\theta_i', 1)$, for $l=1, 2$, with
$\theta_i'=\theta_i/\sqrt{2}$.
We define the estimate of
$T(\theta') = {{k_n}^{-1} \sum_{i=1}^n} |\theta_i'|$ by
%
%
\begin{equation}
\label{sparse.seq.est3}
\widehat{T(\theta')} = {1\over k_n} \sum_{i=1}^n
\bigl\{\tilde\delta(x_{i1}) I\bigl(|x_{i2}|\le2\sqrt{2\log n}\bigr)
+ |x_{i1}| I\bigl(|x_{i2}| > 2\sqrt{2\log n}\bigr)\bigr\},\hspace*{-25pt}
\end{equation}
where $\tilde\delta(\cdot)$ is defined in (\ref{sparse.SK}) and
(\ref{sparse.delta}), and
set the estimator of $T(\theta) = {k_n}^{-1} \times\sum_{i=1}^n |\theta
_i|$ as
%
%
\begin{equation}
\label{sparse.seq.est4}
\widehat{T(\theta)} = \sqrt{2}\widehat{T(\theta')}.
\end{equation}
It can be shown that the estimator $\widehat{T(\theta)}$ is rate
optimal for estimating $T(\theta)$ over~$\Theta_{k_n}$. The proof is
similar to that of Theorem \ref{unbounded.thm}, and
we omit the details here.

\section{Discussions}
\label{discussion.sec}

The present paper was partly inspired by the general theory of
estimating functionals based on i.i.d. observations given in \citet
{DL2} which showed that bounds on minimax estimation
can be based on testing two composite hypotheses.
The difficulty of the composite testing problem was shown in
Le Cam (\citeyear{LeCam73}, \citeyear{LeCam86}) to depend on the total
variation distance between
the convex hulls of the two composite hypotheses.
In the present context the priors $\mu_0$ and $\mu_1$
used in the general lower bound of Section \ref{lowerbd.sec} can be
viewed as
picking two points in the convex hull of two subsets of the parameter
space and Theorem~\ref{lower.bd.thm} gives bounds on the risk over these
two points.
Sections \ref{normal.lowerbd.sec} and \ref{bounded.upperbnd.sec} show
that a careful choice of these priors yields
sharp minimax lower bounds for estimating the $\ell_1$ norm of the
means of
normal random variables.

Best polynomial approximation played a major role in the development
of our results, both for the upper and lower bounds.
Note that the last two conditions in Lemma \ref{prior.lem}
yield
\[
\int\bigl(|t|-G_k^*(t)\bigr) \nu_1(dt) - \int\bigl(|t| - G_k^*(t)\bigr) \nu_0(dt)
= 2 \delta_{k}.
\]
From the definition of $G_k^*$ we have $-\delta_k \le|t| - G_k^*(t)
\le
\delta_k$ for all $-1 \le t \le1$. Since $\nu_0$ and
$\nu_1$
are probability measures it follows that they are
supported on the subsets $A_0$ and $A_1$ of the alternation points
defined in
(\ref{A0}) and (\ref{A1}), respectively.

We should also emphasize that the values of the functional $T$ on the
two sets of support points are not well separated. In fact the values
alternate. This is quite different from the more standard cases of
estimating a linear or quadratic functional. In the case of quadratic
functionals, even though the alternative hypothesis may need to be
composite, the functional only takes on two values, one on the null and
the other on the alternative. See, for example, \citet{CL05}.

The techniques given here can also be compared to those found in
\citet{LNS99} where attention was focused
on estimating the $L_1$ norm of a regression function.
In that paper lower bounds were constructed by mixing in a way
similar to that used in the present paper. However, instead of bounding
a chi-square distance, a bound was given for the Kullback--Leibler
distance.
It is, however, not easy to provide good bounds directly for the
Kullback--Leibler distance. This is particularly true in cases
which correspond to parameter spaces with growing bounds.
The lower bounds provided there only work in the case where
the parameter space has a fixed bound.

For upper bounds, \citet{LNS99}
used a Fourier series approximation of $|x|$,
and the estimate is based on unbiased estimates of
individual terms in the approximation. The maximum error of
the best $K$-term Fourier series approximation can be shown easily to
be of order $K^{-1}$, which is
comparable to the best polynomial approximation of degree $K$.
However, the variance bound of the estimator based on the $K$-term Fourier
series approximation is of order $e^{CK^2}$ for some constant $C>0$,
whereas the variance of our estimator based on the polynomial
approximation of degree~$K$ grows at the rate of $K^K = e^{K\log K}$.
So the variance of the polynomial-based estimator is much smaller
than that of the corresponding estimator using Fourier series, even
though the biases of the two estimators are very similar.
This allows for more terms to be used in the polynomial approximation
with the same variance level thus reducing the bias of the estimate.
In the bounded case, the best rate of convergence for estimators using
Fourier series approximation can be shown to be $(\log n)^{-1}$, which
is sub-optimal relative to the minimax rate $({\log\log n\over
\log n})^2$. Another drawback of
the Fourier series method is that it cannot be used for the unbounded case.

The techniques and results developed in the present paper can be used
to solve other related problems. For example, when the approach taken
in this paper is used for estimating the $L_1$ norm of a regression
function, both the upper and lower bounds given in \citet{LNS99}
are improved. For reasons of space, we shall
report the results elsewhere.
The techniques can also be used for estimating other nonsmooth
functionals such as excess mass. See \citet{CL10}.

\section{Proofs}
\label{proofs.sec}

In this section we first prove the technical lemmas given in the
earlier sections. We then prove Theorem \ref{seq.L1.thm2} in Section
\ref{growing.proof.sec}. The proof of Theorem \ref{unbounded.thm} is
involved and will be given in Section \ref{unbounded.proof.sec}.

\subsection{Proof of technical lemmas}

\mbox{}

\begin{pf*}{Proof of Lemma \ref{prior.lem}}
The proof of this lemma relies on the Hahn--Banach theorem and the
Riesz representation theorem. The argument is essentially the same as
the one given in \citet{LNS99}. We include it here for completeness.

Consider the space $C(-1,1)$ of continuous real-valued functions on the
interval $[-1,1]$ with uniform norm $\| \cdot\|_{\infty}$. Clearly
$f(t) = |t|$ defined on this interval $[-1,1]$ belongs to $C(-1,1)$.
Let $\delta_k$ be the distance in uniform norm on $[-1,1]$ from the function
$f$ to the space of polynomials of order $k$.
Let $\mathcal{P}_k$ be the linear space spanned by the collection of polynomial
of order $k$ and in addition
let $\mathcal{F}_k$ be the linear space spanned by $\mathcal{P}_k$ and
$f$.
Note that every element $g \in\mathcal{F}_k$ can be written uniquely as
$g = cf
+ p_k$ where $p_k \in\mathcal{P}_k$ and $c\in\RR$.
Let $T$ be the linear functional defined by $T(g) = T(cf + p_k) =
c\delta_k$.
It is then clear that $T=0$
on $\mathcal{P}_k$ and $T(f) = \delta_k$.
Now the norm of the functional $T$ is given by
\[
\|T\| \equiv\sup\{ T(g)\dvtx g \in\mathcal{F}_k, \|g\|_{\infty} \le1\}.
\]
It can be checked directly that the norm of this functional is equal to $1$.
Let $G^*_k$ be the closest polynomial in $\mathcal{P}_k$ to $f$.
Then $\|f-G_k^*\|_{\infty} = \delta_k$,
and it follows that ${1 \over\delta_k}(f-G_k^*)$
has a norm of $1$.
Since $T({1 \over\delta_k}(f-G_k^*)) = 1$ it follows that $\|T\| \ge1$.
Now suppose that $\|T\| >1$. Then there exists an
element
$g = cf + p_k$ with $p_k \in\mathcal{P}_k$ such that $\|g\|_{\infty}
=1$ and $T(g)
>1$. This implies that $c>{1 \over\delta_k}$ and
\[
\biggl\| f - \biggl(-{1 \over c}p_k\biggr)\biggr\|_{\infty} = {1 \over c} < \delta_k.
\]
Since $-{1 \over c}p_k\in\mathcal{P}_k$, this is a contradiction to
the definition of $\delta_k$ which is the
distance between $f$ and $\mathcal{P}_k$.

Now by the Hahn--Banach theorem the linear functional $T$ can be
extended to $C(-1,1)$ without increasing the norm of the functional. For
simplicity we shall also call this linear functional $T$.
It then follows from the Riesz representation theorem that for each $g
\in C(-1,1)$
\[
T(g) = \int_{-1}^1 g(t) \tau(d t),
\]
where $\tau$ is a Borel signed measure with total variation equal to $1$.

It follows from Hahn--Jordan decomposition that there exist two positive
measures $\tau_+$ and $\tau_-$ such that $\tau= \tau_+ -\tau_-$.
It then follows that
%
%
\begin{eqnarray}
\label{vanishing0}
\int_{-1}^{1} |t| [ \tau_+(d t) - \tau_-(d t)] &=& \delta_k \quad\mbox
{and}\nonumber\\[-8pt]\\[-8pt]
\int_{-1}^{1} t^l \tau_+(dt) &=& \int_{-1}^{1} t^l \tau_-(dt)\qquad
\mbox{for } l=0,1, \ldots, k.\nonumber
\end{eqnarray}
Define the measures $\tau_-^*$ and $\tau_+^*$ by $\tau_-^*(S) =
\tau_-(-S)$ and $\tau_+^*(S) = \tau_+(-S)$ for all measurable sets $S$.
Then (\ref{vanishing0}) holds with $\tau_-$ and
$\tau_+$ replaced by $\tau_-^*$ and $\tau_+^*$, respectively.
Hence (\ref{vanishing0}) is also true with $\tau_-$ and
$\tau_+$ replaced by $(\tau_- + \tau_-^*)/2$ and $(\tau_+
+\tau_+^*)/2$, respectively. We can thus assume that $\tau$ is symmetric.

Now take $\nu= 2\tau$. Then $\nu$ is symmetric and
%
%
\begin{equation}
\label{vanishing}
\int_{-1}^{1} |t| \nu(dt) = 2\delta_k \quad\mbox{and}\quad
\int_{-1}^{1} t^l \nu(dt) = 0 \qquad\mbox{for } l=0,1, \ldots, k.
\end{equation}

Now let $\nu_{1}$ and $\nu_{0}$ be the positive and the negative
components of $\nu$. Then both $\nu_{1}$ and $\nu_{0}$ are symmetric.
Since $\nu$ has variation equal to $2$ and $\int_{-1}^{1} \nu(dt) = 0$
it follows that $\nu_{1}$ and $\nu_{0}$ are both probability measures.

These measures also clearly satisfy by construction
\[
\int_{-1}^{1} t^l \nu_1(dt) = \int_{-1}^{1} t^l \nu_0(dt)
\]
for $l=0,1, \ldots, k$ and also
\[
\int_{-1}^{1} |t| \nu_1(dt) - \int_{-1}^{1} |t| \nu_0(dt)= 2\delta
_k.
\]
\upqed
\end{pf*}
\begin{pf*}{Proof of Lemma \ref{poly.lem}}
The Chebyshev polynomial $T_{2m}$ can be alternatively written as
%
%
\begin{equation}
T_{2m}(x) = \sum_{l=0}^{m} \Biggl[ (-1)^{m-l} \sum_{j=m-l}^{m}
\pmatrix{2m\cr2j} \pmatrix{j \cr m-l}\Biggr]x^{2l}.
\end{equation}
Write $T_{2m}(x) = \sum_{l=0}^{m} t_{2l}x^{2l}$. Then
%
%
\begin{equation}
\label{TP.bound}
|t_{2l}|= \sum_{j=m-l}^{m} \pmatrix{2m\cr2j} \pmatrix{j \cr m-l}
\le\sum_{j=m-l}^{m} \pmatrix{2m\cr2j} \pmatrix{m \cr m-l} \le2^{2m}2^m
= 2^{3m}.\hspace*{-25pt}
\end{equation}
It is now easy to see that the coefficient for $x^{2k}$ in the
polynomial $G_K(x)$ is bounded from
above by
\[
|g_{2k}| \le{4 \over\pi} \sum_{j=k}^K {2^{3j} \over4j^2 -1} \le
2^{3K}.
\]
The bound on the coefficients $g_{2k}^*$ of the best polynomial
approximation $G_K^*$ follows from Theorem E in \citet{QR07} and
the bound (\ref{TP.bound}).
\end{pf*}
\begin{pf*}{Proof of Lemma \ref{mean-var.lem1}} Write $X = \mu
+ z$ with $z\sim N(0, 1)$. It is well known that $E(H_k^2(z))=k!$,
$E(H_i(z) H_j(z))=0$ for $i\neq j$, and
\[
H_k(\mu+z) = \sum_{j=0}^k \pmatrix{k\cr j} \mu^j H_{k-j}(z).
\]
Hence,
\begin{eqnarray*}
E H_k^2(X)&=&E H^2(\mu+z) = \sum_{i=0}^k \sum_{j=0}^k
\pmatrix{k\cr i}\pmatrix{k\cr j} \mu^{i+j} E(H_{k-i}(z) H_{k-j}(z))\\
&=&\sum_{j=0}^k \pmatrix{k\cr j}^2 \mu^{2j} (k-j)!\\
&=& k! \sum_{j=0}^k \pmatrix{k\cr j} \mu^{2j} {1\over j!}.
\end{eqnarray*}
Note that $k!/j! \le k^{k-j}$ and hence,
\[
E H_k^2(X)= k! \sum_{j=0}^k \pmatrix{k\cr j} \mu^{2j} {1\over j!}
\le k^k \sum_{j=0}^k \pmatrix{k\cr j} \biggl({\mu^{2}\over k}\biggr)^j
= k^k\biggl(1+{\mu^2\over k}\biggr)^k \le e^{\mu^2} k^k.
\]
%
If $|\mu|\le M$ and $M^2 \ge k$, for all $0\le j \le k$,
$\mu^{2j} {1\over j!}\le M^{2j} {1\over j!}\le M^{2k} {1\over k!}$. Hence,
\[
E H_k^2(X)= k! \sum_{j=0}^k \pmatrix{k\cr j} \mu^{2j} {1\over j!}\le
k! \sum_{j=0}^k \pmatrix{k\cr j} M^{2k} {1\over k!} =(2 M^2)^k.
\]
\upqed
\end{pf*}

\subsection{\texorpdfstring{Proof of Theorem \protect\ref{seq.L1.thm2}}{Proof of Theorem 5}}
\label{growing.proof.sec}

For $\theta=(\theta_1,\ldots, \theta_n) \in\RR^n$, denote
\[
b_k(\theta)\equiv{1\over n} \sum_{i=1}^n \theta_i^k.
\]
Note that $E \bar B_k = b_k(\theta)$ for $k \ge0$, and hence
\[
E \widehat{T(\theta)} = \sum_{k=0}^K \tilde g_{2k} b_{2k}(\theta)
={1\over n}\sum_{i=1}^n \tilde G_K(\theta_i).
\]
The bias of
$\widehat{T(\theta)}$ can then be bounded easily as follows:
\begin{eqnarray*}
|E \widehat{T(\theta)} - T(\theta)| &=&\Biggl\vert{1\over n}\sum
_{i=1}^n \tilde
G_K(\theta_i)-{1\over n}\sum_{i=1}^n |\theta_i|\Biggr\vert\le
{1\over n} \sum_{i=1}^n \bigl\vert G_K(\theta_i)- |\theta_i|
\bigr\vert\\
&\le&{2 M_n\over\pi(2K+1)}.
\end{eqnarray*}

Now we consider the variance of $\widehat{T(\theta)}$.
Note that $M_n^2 \ge K$. In this case,
the variance of $\bar B_k$ can be bounded by
\[
\operatorname{Var}(\bar B_{2k}) = n^{-2} \sum_{i=1}^n \operatorname
{Var}(H_{2k}(y_i))
\le n^{-1} (2M_n^2)^{2k}.
\]
Hence
\begin{eqnarray*}
\operatorname{Var}(\widehat{T(\theta)}) &\le&
\Biggl\{\sum_{k=1}^K |\tilde g_{2k}| \operatorname{Var}^{1/2}(\bar B_{2k})
\Biggr\}^2
\le\Biggl\{\sum_{k=1}^K |g_{2k}| M_n^{-2k+1} 2^{k} M_n^{2k}\Biggr\}
^2\cdot n^{-1}\\
&\le& 4 M_n^2 2^{7K} \cdot n^{-1}.
\end{eqnarray*}
With $K = {1\over7}\log_2 n - (\log n)^{1/2}$, the mean squared error
is then bounded by
\[
E\bigl(\widehat{T(\theta)} - T(\theta)\bigr)^2 \le{4 M_n^2\over\pi^2(2K+1)^2}
+ 4 M_n^2 2^{7K} \cdot n^{-1}
= {49 c\over\pi^2} (\log n)^{-1} \bigl(1+o(1)\bigr).
\]

\subsection{\texorpdfstring{Proof of Theorem \protect\ref{unbounded.thm}}{Proof of Theorem 6}}
\label{unbounded.proof.sec}

We now analyze the properties of the hybrid estimator defined in
(\ref{seq.est4}). The key is to study the bias and variance of a~single component.
Let $x_1, x_2\stackrel{\mathrm{i.i.d.}}{\sim} N(\mu, 1)$, and
let
%
%
\begin{equation}\quad
\label{xi}
\xi= \xi(x_1, x_2)=\delta(x_{1}) I\bigl(|x_{2}|\le2\sqrt{2\log n}\bigr)
+ |x_{1}| I\bigl(|x_{2}| > 2\sqrt{2\log n}\bigr).
\end{equation}

Note that
\[
E(\xi) = E\delta(x_{1}) P\bigl(|x_{2}|\le2\sqrt{2\log n}\bigr)
+ E|x_{1}| P\bigl(|x_{2}| > 2\sqrt{2\log n}\bigr).
\]
\begin{lemma}
\label{var.lem1}
Suppose $I(A)$ is an indicator random variable independent~of $X$ and
$Y$, then
%
%
\begin{eqnarray}
\V\bigl(XI(A) + Y I(A^c)\bigr) &=& \V(X) P(A) + \V(Y) P(A^c) \nonumber\\[-8pt]\\[-8pt]
&&{}+ (EX - EY)^2 P(A)P(A^c).\nonumber
\end{eqnarray}
\end{lemma}

Applying Lemma \ref{var.lem1}, we have
%
%
\begin{eqnarray}\label{var.ind}
\V(\xi) &=& \V(\delta(x_{1})) P\bigl(|x_{2}|\le2\sqrt{2\log n}\bigr)
+ \V(|x_{1}|) P\bigl(|x_{2}| > 2\sqrt{2\log
n}\bigr)\hspace*{-20pt}\nonumber\\[-8pt]\\[-8pt]
&&{} +
\bigl(E\delta(x_{1})- E|x_{1}|\bigr)^2 P\bigl(|x_{2}|\le2\sqrt{2\log n}\bigr)
P\bigl(|x_{2}| > 2\sqrt{2\log n}\bigr).\hspace*{-20pt}
\nonumber
\end{eqnarray}

We also need the following lemma for the variance of $\delta$. [The
proof is similar to Lemma 2 in \citet{CL05}.]
\begin{lemma}
\label{var.lem2}
For any two random variables $X$ and $Y$
%
%
\begin{equation}
\V(\min\{X, Y \})
\le
\V X + \V Y.
\end{equation}
In particular, for any random variable $X$ and any constant $C$
%
%
\begin{equation}
\V(\min(X, C)) \le\V X.
\end{equation}
\end{lemma}
\begin{pf}
Without loss of generality we can assume $E(X) = 0$
and $E(Y) \le0$.
Let $Z = \min\{X, Y \}$.
Then
%
%
\begin{equation}
EZ^2 \le EX^2 + EY^2
\end{equation}
and
%
%
\begin{equation}
EZ \le E(Y).
\end{equation}
Hence $(EZ)^2 \ge(EY)^2$ and consequently
%
%
\begin{equation}
\V Z
=
E Z^2 - (EZ)^2
\le EX^2 + EY^2 - (EY)^2
= \V X + \V Y.
\end{equation}
\upqed
\end{pf}
\begin{lemma}
\label{mean-var.lem2}
Let $X\sim N(\mu, 1)$ and $S_K(x) = \sum_{k=0}^K g_{2k} M_n^{-2k+1}
H_{2k}(x)$ with $M_n=8\sqrt{\log n}$ and $K = {1\over12} \log n$.
Then for all $|\mu|\le4\sqrt{2\log n}$,
%
%
\begin{eqnarray}
\bigl|E S_{K}(X) -|\mu|\bigr|&\le& {2 M_n \over\pi(2K+1)},\\
ES_K^2(X) &\le& n^{1/2}\log^5 n.
\end{eqnarray}
\end{lemma}
\begin{pf}
The first part follows from Lemmas \ref{poly.lem} and
\ref{mean-var.lem1} and the discussions in Section \ref{growing.sec}.
To bound $ES_K^2(X)$, it follows from inequality (\ref{inq}) and
Lemmas \ref{poly.lem} and \ref{mean-var.lem1} that
\begin{eqnarray*}
ES_K^2(X) &\le&
\Biggl(\sum_{k=1}^K |g_{2k}| M_n^{-2k+1} (EH_{2k}^2(X))^{1/2}\Biggr)^2
\\
&\le&2^{6K} \Biggl(\sum_{k=1}^K \bigl(8\sqrt{\log n}\bigr)^{-2k+1}
(64\log n)^k\Biggr)^2\\
&\le& n^{1/2}\log^5 n .
\end{eqnarray*}
\upqed
\end{pf}

Write $B(\xi) = E(\xi) - |\mu|$ for the bias of $\xi$.
We divide into three cases according to the value of $|\mu|$. In the
first case when $|\mu|\le\sqrt{2\log n}$, we shall show that the estimator
behaves essentially like $\delta(x_1)$ which is a good estimator when
$|\mu|$ is small. In the second case when $\sqrt{2\log n} \le|\mu|
\le4 \sqrt{2\log n}$, we show that the hybrid estimator uses either
$\delta(x_1)$ or $|x_1|$ and in this case both are good estimators of
$|\mu|$.
In the third case when $|\mu|$ is large, the hybrid estimator is
essentially the same as $|x_1|$.

\textit{Case} 1. \textit{$|\mu|\le\sqrt{2\log n}$.} Note that
$\delta(x_1)$ can be written as $\delta(x_1) = S_K(x_1) -
(S_K(x_1)-n)I(S_K(x_1)\ge n)$ and consequently
%
%
\begin{eqnarray}\label{case1}
|B(\xi)|&=&\bigl|\bigl(E\delta(x_1) P\bigl(|x_{2}| \le2\sqrt{2\log n}\bigr) +
E|x_{1}|\bigr)P\bigl(|x_{2}| > 2\sqrt{2\log n}\bigr) -|\mu|\bigr| \hspace*{-28pt}\nonumber\\
&\le& \bigl|ES_K(x_1)-|\mu|\bigr| + E\bigl\{\bigl(S_K(x_1)-n\bigr)I\bigl(S_K(x_1)\ge n\bigr)\bigr\}
\nonumber\\[-8pt]\\[-8pt]
&&{}+
\bigl(|ES_K(x_1)|+E|x_{1}|\bigr)P\bigl(|x_{2}| > 2\sqrt{2\log n}\bigr) \nonumber\\
&\equiv& B_1 + B_2 + B_3.
\nonumber
\end{eqnarray}
Lemma \ref{mean-var.lem2} yields that
\[
B_1=\bigl|ES_K(x_1)-|\mu|\bigr| \le{2 M_n\over\pi(2K+1)}.
\]
It follows from the fact $|\mu|\le\sqrt{2\log n}$ and the standard bound
for normal tail probability $\Phi(-z) \le z^{-1}\phi(z)$ for $z>0$ that
%
%
\begin{equation}
\label{normal.tail1}
P\bigl(|x_{2}| > 2\sqrt{2\log n}\bigr)\le2 \Phi\bigl(-\sqrt{2\log n}\bigr)\le
{1\over\sqrt{\pi\log n}} n^{-1}.
\end{equation}
Note that in this case
%
%
\begin{eqnarray}
\label{exp1}\quad
|ES_K(x_1)| &=& |\tilde G_K(\mu)| \le|\mu| + {2 M_n\over\pi(2K+1)},
\\
\label{exp2}
E|x_1| &=& |\mu| + 2\phi(\mu)-2 |\mu| \Phi(-|\mu|)\le|\mu| +
1 \le\sqrt{2\log n} + 1.
\end{eqnarray}
It then follows from (\ref{normal.tail1})--(\ref{exp2}) that
\[
B_3\le\biggl(2\sqrt{2\log n} + {2 M_n\over\pi(2K+1)} + 1\biggr)\cdot{1\over
\sqrt{\pi\log n}}
n^{-1}\le3 n^{-1}.
\]
Now consider $B_2$. Note that for any random variable $X$ and any
constant $\lambda> 0$,
%
%
\begin{equation}
E\bigl(XI(X\ge\lambda)\bigr) \le\lambda^{-1} E\bigl(X^2I(X\ge\lambda)\bigr)\le
\lambda^{-1} EX^2.
\end{equation}
This together with Lemma \ref{mean-var.lem2} yields that
%
%
\begin{equation}
\label{B2}
B_2\le E\bigl\{S_K(x_1)I\bigl(S_K(x_1)\ge n\bigr)\bigr\} \le n^{-1} ES_K^2(x_1)
\le n^{-1/2} \log^5 n .
\end{equation}
Combining the three terms together shows that in this case the bias is
bounded by
\[
|B(\xi)| \le B_1 + B_2 + B_3 \le{M_n \over\pi K} \bigl(1+o(1)\bigr).
\]

We now consider the variance. It follows from (\ref{var.ind}) and
Lemma \ref{var.lem2} that
\begin{eqnarray*}
\V(\xi) &\le& \V(S_K(x_{1}))
+ \V(|x_{1}|) P\bigl(|x_{2}| > 2\sqrt{2\log n}\bigr)\\
&&{}+ \bigl(E\delta(x_{1})- E|x_{1}|\bigr)^2 P\bigl(|x_{2}| > 2\sqrt{2\log n}\bigr)\\
&\le& E S_K^2(x_{1}) + Ex_{1}^2 P\bigl(|x_{2}| > 2\sqrt{2\log n}\bigr).
\end{eqnarray*}
Lemma \ref{mean-var.lem2} and equation (\ref{normal.tail1})
together yield that
\[
\V(\xi) \le n^{1/2} \log^5 n \bigl(1+o(1)\bigr).
\]

\textit{Case} 2. $\sqrt{2\log n} \le|\mu| \le4 \sqrt{2\log n}$.
In this case,
\begin{eqnarray*}
|B(\xi)|&=&\bigl|\bigl(E\delta(x_1) P\bigl(|x_{2}| \le2\sqrt{2\log n}\bigr)
 +
E|x_{1}|\bigr)P\bigl(|x_{2}| > 2\sqrt{2\log n}\bigr) -|\mu|\bigr|\\
&\le& \bigl|E\delta(x_1)-|\mu|\bigr| + \bigl|E|x_{1}| -|\mu|\bigr|\\
&\le& \bigl|ES_K(x_1)-|\mu|\bigr| + E\bigl\{\bigl(S_K(x_1)-n\bigr)I\bigl(S_K(x_1)\ge n\bigr)\bigr\} + 2\phi
(\mu).
\end{eqnarray*}
Note that $|ES_K(x_1)-|\mu|| \le{2 M_n\over\pi(2K+1)}$ and as in
(\ref{B2})
\[
E\bigl\{\bigl(S_K(x_1)-n\bigr)I\bigl(S_K(x_1)\ge n\bigr)\bigr\}
\le n^{-1/2} \log^5 n.
\]
Note that $\phi(\mu) \le\phi(\sqrt{2\log n})\le n^{-1}$. Hence,
again the bias is bounded by
\[
|B(\xi)| \le{M_n \over\pi K} \bigl(1+o(1)\bigr).
\]

For the variance, equation (\ref{var.ind}) and
Lemma \ref{var.lem2} yield that
\[
\V(\xi) \le\V(S_K(x_{1})) + \V(x_{1}) + \bigl(E\delta(x_{1})-E|x_{1}|\bigr)^2.
\]
Note that
\begin{eqnarray*}
\bigl(E\delta(x_{1})-E|x_{1}|\bigr)^2 &\le& \bigl[ES_K(x_1) - |\mu| +
E\bigl\{\bigl(S_K(x_1)-n\bigr)I\bigl(S_K(x_1)\ge n\bigr)\bigr\} \\
&&\hspace*{112.5pt}{}- 2\phi(\mu)+2 |\mu|
\Phi(-|\mu|)\bigr]^2 \\
&\le& {M_n^2 \over\pi^2 K^2} \bigl(1+o(1)\bigr).
\end{eqnarray*}
Hence, it follows from Lemma \ref{var.lem2} that
\begin{eqnarray*}
\V(\xi) &\le& ES_K^2(x_{1}) + \V(x_{1}) + \bigl(E\delta
(x_{1})-E|x_{1}|\bigr)^2\\
&\le& n^{1/2} \log^5 n \bigl(1+o(1)\bigr).
\end{eqnarray*}

\textit{Case} 3. $ |\mu| > 4 \sqrt{2\log n}$.
In this case the standard bound
for normal tail probability yields that
\[
P\bigl(|x_{2}| \le2\sqrt{2\log n}\bigr)\le2 \Phi\bigl(-\bigl(|\mu| -2\sqrt{2\log
n}\bigr)\bigr)\le
2 \Phi\biggl(-{|\mu|\over2}\biggr)\le{4\over|\mu|} \phi\biggl({|\mu|\over2}\biggr).
\]
In particular,
\[
P\bigl(|x_{2}| \le2\sqrt{2\log n}\bigr)\le2 \Phi\bigl(-2\sqrt{2\log n}\bigr)\le
{1\over2\sqrt{\pi\log n}} n^{-4}.
\]
Hence,
\begin{eqnarray*}
|B(\xi)|&\le& \bigl|E|x_{1}| -|\mu|\bigr| +
\bigl(|E\delta(x_1)|+E|x_{1}|\bigr)P\bigl(|x_{2}| \le2\sqrt{2\log n}\bigr)\\
&\le& 2\phi(\mu) + (n+ |\mu| + 1)P\bigl(|x_{2}| \le2\sqrt{2\log n}\bigr)\\
&\le& 2\phi(\mu) + 4 \phi\biggl({|\mu|\over2}\biggr) + \hf n^{-3}\le6
\phi\biggl({|\mu|\over2}\biggr) + \hf n^{-3}
\le n^{-3}.
\end{eqnarray*}
For the variance, equation (\ref{var.ind}) and
Lemma \ref{var.lem2} yield that
\begin{eqnarray*}
\V(\xi) &\le& \V(|x_{1}|) + \bigl(\V(\delta(x_{1})) +\bigl(E\delta(x_{1})-
E|x_{1}|\bigr)^2\bigr)P\bigl(|x_{2}|\le2\sqrt{2\log n}\bigr) \\
&\le& 1 + \bigl(3n^2 + 2 (\mu^2 + 1)\bigr)P\bigl(|x_{2}|\le2\sqrt{2\log n}\bigr) =1 + o(1).
\end{eqnarray*}

Putting the three cases together, we have the following.
\begin{proposition}
\label{xi.prop}
For all $\mu\in\RR$,
the bias and the variance of the estimator $\xi$ defined in (\ref{xi})
satisfy
%
%
\begin{equation}
\label{xi.bias.var}
|B(\xi)|\le{M_n\over\pi K}\bigl(1+o(1)\bigr)\quad \mbox{and}\quad
\V(\xi) \le n^{1/2} \log^5 n \bigl(1+o(1)\bigr).
\end{equation}
\end{proposition}
\begin{pf*}{Proof of Theorem \ref{unbounded.thm}} With the detailed analysis
of the one-dimensio\-nal case, we are now ready to give a short proof
Theorem \ref{unbounded.thm}. It suffices to focus on the estimator
$\widehat{T(\theta')}$ given in (\ref{seq.est3}). Note that
\[
\widehat{T(\theta')}= {1\over n} \sum_{i=1}^n \xi(x_{i1}, x_{i2}),
\]
where $\xi$ is defined in (\ref{xi}).
It follows from Proposition \ref{xi.prop} that the bias
$B(\widehat{T(\theta')})$ of the estimator $\widehat{T(\theta')}$ is
bounded by
\[
|B(\widehat{T(\theta')})| \le{1\over n} \sum_{i=1}^n |B(\xi
(x_{i1}, x_{i2}))|
\le{M_n \over\pi K}\bigl(1+o(1)\bigr),
\]
and the variance of $\widehat{T(\theta')}$ is bounded by
\begin{eqnarray*}
\V(\widehat{T(\theta')}) &\le&{1\over n^2} \sum_{i=1}^n \V(\xi
(x_{i1}, x_{i2}))
\le{1\over n^2} \sum_{i=1}^n n^{1/2} \log^5 n \bigl(1+o(1)\bigr)\\
&\le&
64n^{-1/2} \log n \bigl(1+o(1)\bigr).
\end{eqnarray*}
Hence the mean squared error of $\widehat{T(\theta')}$ satisfies
\begin{eqnarray*}
E\bigl(\widehat{T(\theta')} - T(\theta')\bigr)^2 &\le& B^2(\widehat{T(\theta
')}) +
\V(\widehat{T(\theta')}) \le{M_n^2 \over\pi^2 K^2}\bigl(1+o(1)\bigr)\\
&\le&
{C\over\log n} \bigl(1+o(1)\bigr).
\end{eqnarray*}
\upqed
\end{pf*}

\section*{Acknowledgments}

We thank three referees for very constructive comments which have
helped significantly to improve the presentation of the paper.


%
\printaddresses

\end{document}